\documentclass[twocolumn]{autart}

\usepackage{graphicx,
amsmath,amssymb,color}
\usepackage[latin1]{inputenc}
\usepackage[english]{babel}
\usepackage{tikz}

\newcommand{\phitwo}{\phi}

\usepackage{ifthen}
\usepackage{cite}
\usepackage{graphicx}
\usepackage{amsmath}
\usepackage{amssymb}

\newcommand{\beq}[1]{\begin{equation}#1\end{equation}}
\newcommand{\bal}[1]{\begin{align}#1\end{align}}
\newcommand{\bpm}[1]{\begin{pmatrix}#1\end{pmatrix}}

\newcommand{\fci}[3]{\begin{figure}\begin{center} \includegraphics[#1]{{#2}} \caption{{#3}} \end{center}\end{figure}}
\newcommand{\dom}{\ensuremath{\mbox{dom }}}

\newtheorem{dfn}{Definition}

\newcommand{\Rl}[2]{\ensuremath{\mathbb{R}^{#1}_{#2}}}   

    \renewcommand*{\theenumi}{(\roman{enumi})}

\begin{document}
\begin{frontmatter}
\title{On the graphical stability of hybrid solutions with non-matching jump times: Extended Paper}
\thanks{Corresponding author J.~J.~B.~Biemond.}
\author[TNO]{J.~J.~Benjamin~Biemond}\ead{benjamin.biemond@tno.nl},
\author[ULOR]{Romain Postoyan}, 
\author[TUE]{W.~P.~Maurice~H.~Heemels}, 
\author[TUE,UoM]{Nathan van de Wouw}

\address[TNO]{Department of Optomechatronics, Netherlands Organization for Applied Scientific Research, TNO, Delft, the Netherlands.}
\address[ULOR]{Universit\'e de Lorraine, CNRS, CRAN, F-54000 Nancy, France}
\address[TUE]{Department of Mechanical Engineering, Eindhoven University of Technology, the Netherlands}
\address[UoM]{Department of Civil, Environmental \& Geo- Engineering, University of Minnesota,
Pillsbury Drive SE, U.S.A.}

\begin{abstract}
We investigate stability of a solution of a hybrid system in the sense that the graphs of solutions from nearby initial conditions remain close and tend towards the graph of the given solution. In this manner, a small continuous-time mismatch is allowed between the jump times of neighbouring solutions and the `peaking phenomenon' is avoided. We provide conditions such that this stability notion is implied by stability with respect to a specifically designed distance-like function. Hence, stability of solutions in the graphical sense can be analysed with existing Lyapunov techniques.
\end{abstract}

\end{frontmatter}

\section{Introduction}\label{sect-intro}
Hybrid systems feature both continuous evolution in time and discrete events, and are valuable for the modelling, analysis and control of many engineering applications, see \cite{goe_san_12,lei_wouw_08sh} and the references therein. While the stability of stationary points and sets for hybrid systems is relatively well understood, far less is known about the stability of a given time-varying and jumping solution. 
The stability of time-varying solutions to hybrid systems is challenging \cite{lei_wouw_08sh,mor_bro_10} as two nearby solutions typically show `peaking behaviour', i.e. they experience jumps at close, but not identical, jump times and during this time-mismatch interval, the state distance between both solutions will generally not be small. To analyse stability of a hybrid solution, set-stability techniques \cite{for_teel_13sh}, ignoring the state difference in this interval \cite{mor_bro_10}, or non-Euclidean distance-like functions \cite{bie_hee_16,bie_wouw_13sh} have been proposed. The approaches in \cite{for_teel_13sh,mor_bro_10} seem to be hard to generalize, while stability with respect to non-Euclidean distance-like functions is hard to interpret.
Here, we investigate the stability of a \emph{single} solution (in contrast to incremental stability \cite{li_phi_16,bie_pos_18sh}) using closeness of the graphs of solutions. We define stability of a given solution \emph{in graphical sense}, which implies that the Euclidean distance between the states of this solution and nearby solutions, when compared at close continuous-time instances, tends towards zero and, in addition, the difference between the continuous times used in this comparison also tends to zero when time evolves (a similar definition is used in \cite{san_teel_07}). 
We prove that graphical stability is implied by stability with respect to a well-designed non-Euclidean distance-like function, cf.\ \cite{bie_hee_16,bie_wouw_13sh}. This implication allows to prove (asymptotic) stability in graphical sense using existing Lyapunov-based techniques.
For the first time, a graphical and intuitive notion of stability is provided for a given solution that applies to a large class of hybrid systems and can be analysed with existing techniques. 

\section{Hybrid system and stability definitions}
Let $(x,y)=[x^{\mathrm{T}},\,y^{\mathrm{T}}]^{\mathrm{T}}$ for $(x,y)\in\Rl{n}{}\times \Rl{m}{}$ and for
 a set-valued mapping $F$, $\dom F := \{x\in\Rl{n}{}\,:\,F(x)\neq\emptyset\}$. Given $S_1,S_2\subset \mathbb{R}^n$, $S_1+S_2$ denotes $\{y_1+y_2: y_1\in S_1,y_2\in S_2\}$ and  $\mathcal{B}_s$ denotes $\{x\in \mathbb{R}^n\,:\,\|x\|\leq s\}$, with $\|\cdot \|$ the Euclidean norm.
Let $\rho_{B_1}(x,y):=\inf_{(u,v)\in {B_1}}\|(x-u,y-v)\|$ for $B_1\subset \mathbb{R}^{n}\times \mathbb{R}^n$.
We study hybrid systems
\begin{equation}
\dot{x}  \in  F(x),  x \in C; \quad x^{+}  \in G(x), \in D,
\label{eqhs}
\end{equation}
as in \cite{goe_san_12} where we impose the \emph{hybrid basic conditions} (cf.\ \cite{goe_san_12}) given by:\\
(A1) flow set $C$ and jump set $D$ are closed subsets of $\Rl{n}{}$\\
(A2) flow map $F:\Rl{n}{}\rightrightarrows\Rl{n}{}$ is outer semicontinuous 
and locally bounded 
relative to $C$, $C\subset \dom F$, and $F(x)$ is convex for each $x\in C$\\
(A3) jump map $G:\Rl{n}{}\rightrightarrows\Rl{n}{}$ is outer semicontinuous and locally bounded relative to $D$, and $D\subset \dom G$.\\
Let hybrid time domain, maximal solutions to \eqref{eqhs} and tangent cone at a point $x\in \mathbb{R}^n$ to a set $B\subseteq\mathbb{R}^n$ (denoted $T_B(x)$) be defined as in \cite{goe_san_12}.
We call solution  $\phi$ $t$-complete if $\sup_{(t,j)\in \dom \phi}t=\infty$ and bounded when $\|\phi(t,j)\|\leq R$ for all $(t,j)\in \dom \phi$ and some $R>0$.\\
%

Exploiting the Hausdorff distance between solution graphs (see \cite{goe_teel_06,san_teel_07}), we define graphical stability of a solution as follows: 
\begin{dfn}\label{defstabgraf} A $t$-complete solution $\phi^\star$ to system \eqref{eqhs} is \emph{stable in graphical sense} if the following condition holds. For any $\varepsilon>0$, there exists $\delta(\epsilon)>0$ such that for any maximal solution $\phitwo$ with $\|\phi^\star(0,0)-\phitwo(0,0)\|< \delta(\epsilon)$ it holds that\\
(i) for all $(t,j)\in\dom\phi^\star$, there exists $(t',j')\in\dom\phitwo$ with $|t-t'|<\varepsilon$ such that $\|\phi^\star(t,j)-\phitwo(t',j')\|<\varepsilon$,\\
(ii) for all $(t',j')\in\dom\phitwo$, there exists $(t,j)\in\dom\phi^\star$ with $|t-t'|<\varepsilon$ such that $\|\phi^\star(t,j)-\phitwo(t',j')\|\!<\!\varepsilon$.\\
The solution $\phi^\star$ is \emph{asymptotically stable in graphical sense} if, in addition, there exists $r>0$ such that for any $\varepsilon>0$ and any maximal solution $\phitwo$ with $\|\phi^\star(0,0)\!-\!\phitwo(0,0)\|\!<\!r$ there exists $T\!\geq\! 0$ for which the following statements hold:\\
(iii) for all $(t,j)\in\dom\phi^\star$ with $t\geq T$, there exists $(t',j')\in\dom\phitwo$ with  $|t-t'|< \varepsilon$ such that $\|\phi^\star(t,j)-\phitwo(t',j')\|\!<\! \varepsilon$,\\
(iv) for all $(t',j')\in\dom\phitwo$ with $t'\geq T$, there exists $(t,j')\in\dom\phi^\star$ with  $|t-t'|< \varepsilon$ such that $\|\phi^\star(t,j)-\phitwo(t',j')\|\!<\! \varepsilon$.\hfill $\Box$
\end{dfn}
Definition~\ref{defstabgraf} prioritises continuous time over the jump counter (in \cite{bie_pos_18sh} and references therein, this is used for incremental stability). Focussing on hybrid systems that cannot exhibit consecutive jumps without flow, \cite[Definition 1]{bie_hee_16} yields the distance-like function $
\rho_{\mathcal{A}}(x_1,x_2)$  with
$\mathcal{A}:=\big\{(x_1, x_2)\in (C\cup D\cup G(D))^2\,:\,x_1=x_2$ or $x_2\!\in\! D\!\land\! x_1\!\in\! G(x_2)\mbox{\! or\! } x_1\!\in\! D\!\land \!x_2\!\in\! G(x_1)\big\}$ and by  \cite[Lemma 1]{bie_pos_18sh} the stability definition in \cite{bie_wouw_13sh} is  equivalent to:
\begin{dfn}\label{dfnstabtrajrhoa} A $t$-complete solution $\phi^\star$ to system \eqref{eqhs} is \emph{stable with respect to $\rho_{\mathcal{A}}$} if the following conditions hold. For any $\varepsilon_w>0$, there exists $\delta_w(\varepsilon_w)>0$ such that for any maximal solution $\phitwo$ with $\rho_{\mathcal{A}}(\phi^\star(0,0),\phitwo(0,0))< \delta_w(\varepsilon_w)$ it holds that\\
(i) for all $(t,j)\in\dom\phi^\star$, there exists $(t,j')\in\dom\phitwo$  such that $\rho_{\mathcal{A}}(\phi^\star(t,j),\phitwo(t,j'))<\varepsilon_w$, and\\
(ii) for all $(t,j')\in\dom\phitwo$, there exists $(t,j)\in\dom\phi^\star$ such that $\rho_{\mathcal{A}}(\phi^\star(t,j),\phitwo(t,j'))\!<\!\varepsilon_w$.\\
The solution $\phi^\star$ is \emph{asymptotically stable with respect to $\rho_{\mathcal{A}}$} if, in addition,  there exists $r_w\!\!>\!\!0$ such that for any $\varepsilon_w\!\!>\!\!0$ and any maximal solutions $\phitwo$ with $\rho_{\mathcal{A}}(\phi^\star(0,0),\phitwo(0,0))\!\!<\!\! r_w$ there exists $T_w\!\!\geq\!\! 0$ for which it holds that\\
(iii) for all $(t,j)\in\dom\phi^\star$ with $t\geq T_w$, there exists $(t,j')\in\dom\phitwo$ such that $\rho_{\mathcal{A}}(\phi^\star(t,j),\phitwo(t,j'))\!<\! \varepsilon_w$, and\\
(iv) for all $(t,j')\in\dom\phitwo$ with $t'\geq T_w$, there exists $(t,j')\in\dom\phi^\star$ such that $\rho_{\mathcal{A}}(\phi^\star(t,j),\phitwo(t,j'))\!<\! \varepsilon_w$.\hfill $\Box$
\end{dfn}

\section{Comparison of stability concepts\label{secmainthm}}
The following result allows to compare both definitions.
\begin{thm}\label{thm-mismatch}
Consider system \eqref{eqhs} with a $t$-complete solution $\phi^\star$ and let the following conditions hold:\\
(i) $G(D)\cap D=\emptyset$, $G(D)\subset C$ and $G$ is single-valued and proper;\\
(ii) $
    \forall x\in C\cap D,\  \phantom{-}F(x)\cap T_C(x)=\emptyset;$\\
(iii) $\forall x\in C\cap G(D),\ -F(x)\cap T_C(x)=\emptyset$;\\
(iv) either $D$ is bounded or $\phi^\star$ is bounded.\\
Then, for all $\varepsilon>0$ there exists $s>0$ such that for any $t$-complete solution $\phitwo$ to \eqref{eqhs} that satisfies the conditions:\\
(v) $\|\phi^\star(0,0)-\phitwo(0,0)\|< s;$\\
(vi) for all $(t,j)\in \dom \phi^\star$, there exists $ (t,\tilde{\jmath})\in \dom \phitwo$ such that $\rho_{\mathcal A}(\phi^\star(t,j),\phitwo(t,\tilde{\jmath}))< s$,\\
it holds that for any $(t,j)\!\!\in\!\!\dom\phi^\star$, there exists $(t',j')\!\!\in\!\!\dom\phitwo$ with $|t-t'|\!\!<\!\varepsilon$ such that $\|\phi^\star(t,j)-\phitwo(t',j')\|\!\!<\!\varepsilon$.
\hfill $\Box$
\end{thm}
The existence of a solution $\phitwo$ verifying (v) and (vi) is guaranteed by \cite[Proposition 6.14]{goe_san_12} and that the conclusion of Theorem \ref{thm-mismatch} coincides with (i) of Definition~\ref{defstabgraf}.
Condition (i) in this theorem restricts the possibility of forward and backward jumps of the hybrid system (e.g., excluding Zeno-type solutions). With the definition of the set $\mathcal A$, this condition is essential to draw conclusions on the difference between solutions $\phitwo$ and $\phi^\star$ when condition (vi) holds. Condition (ii), combined with the fact that $F(x)$ is nonempty for $x\in C$,  guarantees that solutions that are close to the jump set, will indeed jump in the near future (with a uniform bound on the jump time mismatch). In particular, it ensures $C\cap D$ has zero Lebesgue measure and in case $D$ is a submanifold with $C$ located on one side of this manifold, (ii) enforces transversal intersection of solutions with $D$. Condition (iii) has a similar role for solutions backward in time and also guarantees that solutions cannot enter $G(D)$ by flow. To infer compactness results for a bounded subset of the jump set $D$ including those points of $D$ explored by the solution $\phi^\star$, condition (iv) is imposed.
We now formulate our main result below.
\begin{thm}\label{thm-dafis}
Consider a $t$-complete solution $\phi^\star$ to system \eqref{eqhs}, suppose that conditions (i)-(iv) of Theorem~\ref{thm-mismatch} hold.
If the solution $\phi^\star$ is (asymptotically) stable with respect to $\rho_{\mathcal A}$ (as in Definition~\ref{dfnstabtrajrhoa}), then it is (asymptotically) stable in graphical sense as in Definition~\ref{defstabgraf}.\hfill $\Box$
\end{thm}

%

\section{Example \label{secexample}}
Consider a single degree-of-freedom mechanical system with unit mass, damper with  damping constant $0.02$ and spring with unit stiffness constant and unloaded position $x_1=1$. The dynamics near a reference solution $\phi^\star$ is 
given as
\bal{
\dot x&=\left(\begin{smallmatrix}0& 1\\ -1&-0.02\end{smallmatrix}\right)x+\left(\begin{smallmatrix}0\\k\bar x_1\end{smallmatrix}\right)+\left(\begin{smallmatrix}0\\1\end{smallmatrix}\right)(u_{\text{ff}}(t)+u_{\text{fb}}(t,x)),
\nonumber \\&\hspace{20mm}
x\in C=\{(z_1,z_2)\in \mathbb{R}^2\,:\,z_1\geq 0\},\nonumber \\
x^+&=-\varepsilon x,\hfill
x\in D=\{(z_1,z_2)\in \mathbb{R}^2\,:\,z_1=0,z_2\leq -r\}. \nonumber
}
Let the forcing $u_{\text{ff}}(t)$ and feedback $u_{\text{fb}}(t,x)$ be selected as in \cite[Section 6]{bie_hee_16}, where we note that this solution $\phi^\star$ is bounded and $t$-complete, and using a Lyapunov function argument, in \cite{bie_hee_16} the applied feedback is proven to render the solution $\phi^\star$ asymptotically stable with respect to $\rho_{\mathcal{A}}$. Item (i) of Theorem~\ref{thm-mismatch} holds as $G(D)=\{(z_1,z_2)\in \mathbb{R}^2\,:\,z_1=0,z_2\geq \varepsilon r\}$, and (ii) and (iii) are verified since $T_C(x)=\{z\in \mathbb{R}^2\,:\, \bpm{0&1}z\geq 0\}$ for $x\in D\cup G(D)$ and $\bpm{0&1}F(x)=x_2$. Consequently, Theorem~\ref{thm-dafis} is applicable and ensures that $\phi^\star$ is (asymptotically) stable in graphical sense.
\\
In Fig.~\ref{figeps}, solutions of this hybrid system are shown. The function $\rho_{\mathcal{A}}$ is shown in panel c) and illustrates that the solution $\phi^\star$ is asymptotically stable with respect to $\rho_{\mathcal{A}}$. Indeed, as stated in Theorem~\ref{thm-dafis}, $\phi^\star$ is asymptotically stable in graphical sense, see panels a)-b).
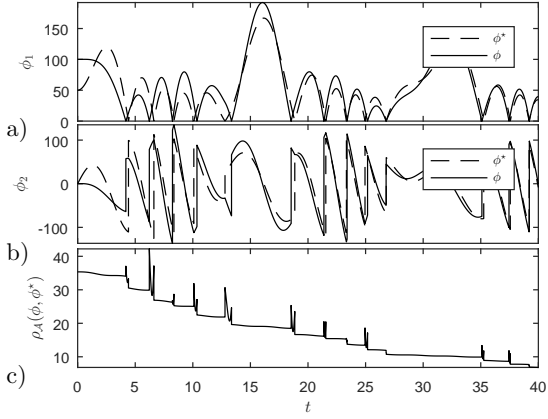
\fci{width=.85\columnwidth}{figeps3}{a) and b) Reference solution $\phi^\star$ and a nearby solution  $\phi$. c)~Distance function $\rho_{\mathcal{A}}$. \label{figeps}}

%
Analysis of stability in graphical sense is facilitated by Theorem~\ref{thm-dafis}, which states that this stability notion is implied by stability of the solution with respect to a specifically constructed distance-like function. Hence, existing Lyapunov-based approaches as in \cite{bie_wouw_13sh} can be used to prove asymptotic stability in graphical sense. The example illustrates how Theorem~\ref{thm-dafis} is used to prove stability in graphical sense for a bouncing ball tracking problem. An open question is if the stability definitions in graphical sense, or in terms of $\rho_{\mathcal{A}}$, are equivalent.

\bibliographystyle{ieeetran}
\bibliography{IEEEabrv,c:/postdoc_files/kuleuven/literature/allrefs}

\appendix
\section{Proofs of Theorems~\ref{thm-mismatch} and \ref{thm-dafis}}
%
\begin{lem}\label{lemprej}
Consider system \eqref{eqhs}, suppose  (ii) of Theorem~\ref{thm-mismatch} holds and let $\bar K>0$ be given.
If $D$ is bounded then for all $\epsilon_2>0$, there exists $\epsilon_1>0$ such that for any $t$-complete solution $\phi$ to \eqref{eqhs} and any $(t,j)\in \dom \phi$ such that $\phi(t,j)\in (C\cap D)+\mathcal{B}_{\epsilon_1}$,
there exists $t'\in [t,t+\epsilon_2]$ such that $
(t',j)\in \dom \phi$ and $\phi(t',j)\in C\cap D.$
For unbounded $D$, such $\epsilon_1$ and $t'$ exist if $\phi(t,j)\in \mathcal{B}_{\bar K}$.$\hfill \Box$
\end{lem}
\begin{pf}
Let $\widehat D$ be given by $D$ if it is bounded and by $D\cap \mathcal B_{2\bar K}$ otherwise. Since $D$ is closed by (A1), $\widehat D$ is compact.
Since $F$ is locally bounded by (A2) of the hybrid basic assumptions, we can select $\gamma>0, \tilde F>1$ such that $\|f\|<\tilde F$ for all $f\in F(x)$ and $x\in (C\cap \widehat D)+\mathcal B_\gamma$.
Given  $\epsilon_2>0$, we show that  $\epsilon_1>0$ can be selected such that for any solution $\phi$ and any hybrid time $(t,j)\in \dom \phi$ such that $\phi(t,j)\in (C\cap D)+\mathcal{B}_{\epsilon_1}$, we have
\begin{align}
\|\phi(s,j)\!-\! \phi(t,j)\|\leq |s\!-\!t|\tilde F\leq \epsilon_2,\ \mbox{for all }(s,j)&\in \dom \phi\nonumber\\
\mbox{ with } s\in [t,t+\min(\tfrac{\gamma}{1+\tilde F},\tfrac{\epsilon_2}{\tilde F})).&\label{eq-lem-phiclosegam}
\end{align}
Namely, impose $0< \epsilon_1<\min(\bar K,\tfrac{\gamma}{1+\tilde F})$ and consider a $t$-complete solution $\phi$ (with $\|\phi(t,j)\|\leq \bar K$ if $D$ is unbounded due to the hypothesis of the lemma) and $(t,j)\in \dom \phi$ with
$\phi(t,j)\in (C\cap D)+\mathcal{B}_{\epsilon_1}$. Exploiting $\epsilon_1<\bar K$ and $\|y\|> 2\bar K$ for all $y\in  D\setminus \widehat D$ in case $D$ is unbounded, we find $\phi(t,j)\in (C\cap \widehat D)+\mathcal{B}_{\epsilon_1}.$
 Introducing $\delta t:=\min(\tfrac{\gamma}{1+\tilde F},\tfrac{\epsilon_2}{\tilde F})$, we let $\bar s\in [t,t+\delta t]$ be the maximum  scalar such that for all $s\in [t,\bar s]$, $(s,j)\in \dom \phi$ and $\phi(s,j) \in \{\phi(t,j)\}+\mathcal B_{\gamma-\epsilon_1}$ hold.
We deduce $\phi(s,j)\in (C\cap \widehat D)+\mathcal{B}_\gamma$ for $s\in [t,\bar s]$ and, hence, $\|f\|< \tilde F$ for all $f\in F(\phi(s,j))$ and $s\in [t,\bar s]$. By definition of the solution, $\frac{d\phi(s,j)}{ds}\in F(\phi(s,j))$ for almost all $s\in [t,\bar s]$, such that
$\|\phi(s,j)-\phi(t,j)\|\leq  | s-t|\tilde F,\mbox{ for }s\in [t,\bar s]$
holds.
This directly implies that for $s\in [t,\bar s]$, $\|\phi(s,j)-\phi(t,j)\|\leq \delta t \tilde F\leq \tfrac{\gamma}{1+\tilde F}\tilde F =\gamma-\tfrac{\gamma}{1+\tilde F}< \gamma-\epsilon_1$. Hence, we find $\bar s=\min\{t+\delta t,\max(s\in \mathbb{R}\,:\, (s,j)\in \dom \phi)\}$ and $[t,\bar s]$ coincides with $\{s\in [t,t+\delta t]\,:\, (s,j)\in \dom \phi\}$. With $s$ as above and $|s-t|\leq \delta t\leq \tfrac{\epsilon_2}{\tilde F}$,
\eqref{eq-lem-phiclosegam} is proven.\\
For the sake of contradiction, we now suppose:
\begin{itemize}
\item[S1:]
  there exists $\epsilon_2>0$ such that for all $\epsilon_1\in (0,\gamma)$, there exists a $t$-complete solution $\phi$ to system \eqref{eqhs}, with $\phi(0,0)\in (C\cap \widehat D)+\mathcal{B}_{\epsilon_1}$, and $\phi(t',0) \notin C\cap D$ for all $t'\in [0,\min(\tfrac{\gamma}{1+\tilde F},\tfrac{\epsilon_2}{\tilde F})]$.
\end{itemize}
Let  $\epsilon_2>0$ be as in S1.
We select $\epsilon_1\in (0,\gamma)$ such that \eqref{eq-lem-phiclosegam} holds for any solution $\phi$ and hybrid time instant $(t,j)\in \dom \phi$ for which $\phi(t,j)\in C\cap \widehat D+\mathcal{B}_{\epsilon_1}$ and introduce $\delta t':=\min(\tfrac{\gamma}{1+\tilde F},\tfrac{\epsilon_2}{\tilde F},\tfrac{\gamma-\epsilon_1}{\tilde F})$. Select $\epsilon_1^0=\epsilon_1$, a sequence $\{\epsilon^i_1\}_{i\in \{0,1,\ldots\} }$ of strictly positive scalars such that $\epsilon_1^i<\epsilon_1^{i-1}$ for all $i\in \{1,2,\ldots\}$, $\lim_{i\to \infty}\epsilon_1^i=0$, and a sequence $\{\phi^i\}_{i\in \{1,2,\ldots\}}$  of solutions to \eqref{eqhs} with $\phi^i:[0,\delta t']\times \{0\}\to \mathbb{R}^n$, $i\in \{1,2,\ldots\}$,
$
\phi^i(0,0)\notin C\cap D$ ;
$\phi^i(0,0)\in (C\cap \widehat D)+\mathcal{B}_{\epsilon_1^i}$ and $\phi^i(t',0)\in (\{\phi^i(0,0)\}+\mathcal B_{\epsilon_2})\setminus (C\cap D)\mbox{ for all }t'\in [0,\delta t']
$
(cf.\ S1).
Introducing the compact set $\mathbb K=(C\cup \widehat D)\cap (C\cap \widehat D+\mathcal{B}_\gamma)$ we find with \eqref{eq-lem-phiclosegam} and $\delta t'\leq \tfrac{\gamma-\epsilon_1^0}{\tilde F}$ that for each solution $\phi^i$, $i\in \mathbb{N}$, it holds that $\phi^i(t',0) \in \mathbb{K}$ for $t'\in[0,\delta t']$, where we exploited the bound $\|\frac{d\phi(t,j)}{dt}\|\leq \bar F$ that holds in this time interval. Hence, the elements of the sequence $\{\phi^i\}_{i\in \mathbb{N}}$ are contained in the bounded set of absolutely continuous functions $[0,\delta t']\times \{0\}\to \mathbb K$. Consequently, there exists a convergent subsequence within $\{\phi^i\}_{i\in \mathbb{N}}$ that graphically converges 
to a function $\phi:[0,\delta t']\times \{0\}\to \mathbb K$ with $\phi(0,0)\in C\cap \widehat D$ according to Theorem~5.7 in \cite{goe_san_12}. This function $\phi$ is a solution to the hybrid system by sequential compactness of solutions to hybrid systems satisfying A1)-A3), cf.~\cite[Theorem 6.8 and Definition 6.2(a)]{goe_san_12}. However, the existence of such a solution is excluded by item (ii) of Theorem~\ref{thm-mismatch} as no solutions to \eqref{eqhs} can flow on $C\cap D$. Hence, S1 is contradicted and we have proven the lemma for $(t,j)=(0,0)$. Time-invariance of \eqref{eqhs} concludes this proof.
$\hfill \Box$
\end{pf}
\begin{lem}\label{lempostj}
Consider a hybrid system \eqref{eqhs} satisfying items (i), (iii) of Theorem~\ref{thm-mismatch} and let $\bar K\in \mathbb{R}$ be given.
If $D$ is bounded then for all $\epsilon_4>0$, there exists $\epsilon_3>0$ such that for any solution $\phi$ to \eqref{eqhs} and any $(t,j)\in \dom \phi$ such that $\phi(t,j)\in (C\cap G(D))+\mathcal{B}_{\epsilon_3}$, $j\in \{1,2,\ldots\}$,
there exists $t'\in [t-\epsilon_4,t]$ such that
$
(t',j)\in \dom \phi$ and $\phi(t',j)\in C\cap G(D).$
For unbounded $D$, such $\epsilon_3$ and $t'$ exist if $\phi(t,j)\in \mathcal{B}_{\bar K}$. $\hfill \Box$
\end{lem}
\begin{pf}  If $j>0$, we observe that the flowing solution segment of $\phi$ to \eqref{eqhs} is characterised by the differential inclusion $x'\in -F(x),x\in C$ as long as $x\not \in G(D)$ and the direction of continuous time is reversed. Hence, we deduce that the statement of Lemma~\ref{lempostj} is proven by application of Lemma~\ref{lemprej} after replacing $D$ with $G(D)$.$\hfill \Box$\end{pf}
%
\def\Elproofname{PROOF of THEOREM~\ref{thm-mismatch}.}
\begin{pf}
If the jump set $D$ is unbounded (cf.\ (iv)) we construct $\bar K>0$ that verifies Lemmas~\ref{lemprej} and \ref{lempostj} and  is such that $\rho_{\mathcal A}(\phi^\star(t,j),\phitwo(t,\tilde{\jmath}))$ in (vi) can be written as the distance from a compact set when both (iv) and (vi) hold for some $s<K$. For this purpose, we select $\bar K>0$ such that both $G^{-1}(G(D)\cap \mathcal{B}_{2K})+\mathcal B_K\subset D\cap \mathcal B_{\bar K}$ and $(G(D\cap \mathcal B_{2K})+\mathcal{B}_K)\subseteq \mathcal B_{\bar K}$ hold.
We define $\widehat D=D\cap \mathcal{B}_{\bar K}$ if $D$ is unbounded and $\widehat D=D$ otherwise. The set $\widehat D$ is closed by (A1), it is thus compact. In addition, $
\widehat G_D:=G(\widehat D)$ is compact by locally boundedness and outer semi-continuity of $G$.
\\
We introduce $\mathcal{A}_{01}=\{(z_1,z_2)\in (C\cup D)^2\,:\, z_2= G(z_1),\ z_1\in \widehat D\}$ such that the set $\mathcal{A}_{01}$ is compact as $G$ is locally bounded and outer semi-continuous. Introducing the symmetrical set $\mathcal{A}_{10}=\{(z_1,z_2)\in (C\cup D)^2\,:\, z_1= G(z_2),\ z_2\in \widehat D\}$, compactness of this set follows from the symmetry. From item (i) of Theorem~\ref{thm-mismatch}, we conclude $\mathcal{A}_{01}$ and $\mathcal{A}_{10}$ are not intersecting. Furthermore, the intersection of $\mathcal{A}_{01}$ and $\mathcal{A}_{00}:=\{(z_1,z_2)\in (C\cup D)^2\,:\, z_2=z_1\}$ is empty, as, for points $(z_1,z_2)$ in this intersection, $z_1=z_2= G(z_1)$ should hold, contradicting item (i); similarly, we find $\mathcal A_{10}\cap \mathcal A_{00}=\emptyset$. Since  $\mathcal{A}_{00}$, $\mathcal A_{01}$, $\mathcal A_{10}$ are disconnected, closed and the latter two sets compact, there exists $\bar s>0$ such that $\{(x,y)\in (C\cup D)^2\,:\, \rho_{\mathcal A_{00}}(x,y)\leq \bar s\}$, $\{(x,y)\in (C\cup D)^2\,:\, \rho_{\mathcal A_{01}}(x,y)\leq \bar s\}$ and $\{(x,y)\in (C\cup D)^2\,:\,  \rho_{\mathcal A_{10}}(x,y)\leq \bar s\}$ are mutually disconnected.
As $\widehat D$ and $\widehat G_D$ are compact and $F$ is locally bounded by (A2), there exist positive scalars $r,\tilde F$ such that $\tilde F>1$ and $\|f\|\leq \tilde F$ for all $f\in F(x)$ and $x\in (\widehat D\cup \widehat G_D)+\mathcal B_r$.
Now, fix $\varepsilon>0$ as in Theorem~\ref{thm-mismatch} and take $\epsilon_1$ as in Lemma~\ref{lemprej} with $\epsilon_2=\frac{\min(\varepsilon,2r)}{2(1+\tilde F)}$ and  $\epsilon_3$ as in Lemma~\ref{lempostj} with $\epsilon_4=\frac{\min(\varepsilon,2r)}{2(1+\tilde F)}$, where, if $D$ is unbounded, $\bar K$ is used.
\\
Selecting $s>0$ such that $
s< \min\big(\bar s,\tfrac{\varepsilon}{2\sqrt 2\tilde F},\tfrac{r}{\tilde F+1},K,$ $\tfrac 12 \min_{{u\in \widehat D,w\in \widehat G_D}}\|u-w\|,\epsilon_1,\epsilon_3\big)$
(with $K=\infty$ if $D$ is bounded), we prove that condition~(i) in Definition~\ref{defstabgraf} holds. Considering any pair $(\phi^\star,\phitwo)$ of $t$-complete solutions to \eqref{eqhs} satisfying (v),(vi) and selecting $(t,j)\in \dom \phi^\star$ arbitrary, we find by (vi) that there exists $(t,\tilde{\jmath})\in \dom \phitwo$ such that $\rho_{\mathcal{A}}(\phi^\star(t,j),\phitwo(t,\tilde{\jmath}))<s$. Exploiting strictness of this inequality and the infimum defining $\rho_{\mathcal{A}}$, this implies that there exists $(z_1,z_2)\in \mathcal{A}$ such that
\beq{\label{eqphis}
\|(\phi^\star(t,j)-z_1,\phitwo(t,\tilde{\jmath})-z_2)\|<s
}
holds, with $z_1,z_2$ satisfying one of the following three cases that are generated by the `or' conditions in the definition of $\mathcal{A}$. We now construct $(t',j')\in \dom \phitwo$ as in the theorem.
\\
\textbf{Case 1:  $z_1=z_2\in C\cup D$.\quad } We directly observe $(z_1,z_2)\in \mathcal{A}_{00}$ and select $(t',j')=(t,\tilde{\jmath})$. From $\|(\phi^\star(t,j)-z_1,\phitwo(t,\tilde{\jmath})-z_2)\|\geq$ \\ $ \displaystyle\min_{z\in \mathbb{R}^n}
\|(\phi^\star(t,j)-z,\phitwo(t,\tilde{\jmath})-z)\|=\tfrac{1}{\sqrt 2} \left\|\phi^\star(t,j)-\phitwo(t,\tilde{\jmath})\right\|$ and \eqref{eqphis}, we conclude $\|\phi^\star(t,j)-\phitwo(t,j')\|<\sqrt 2 s\leq \varepsilon$ since $s<\tfrac \varepsilon{2\sqrt 2\tilde F}< \tfrac {\varepsilon} {\sqrt{2}} $. Hence, $(t',j')$ satisfies the conditions imposed in the theorem.
\\
\textbf{Case 2: $z_1\in D, z_2= G(z_1)$.\quad }
Since $\phitwo(t,\tilde \jmath)$ is close to $G(D)$, we will apply Lemma~\ref{lempostj} to prove that $\phitwo$ experienced a jump shortly before the time instant $(t,\tilde \jmath)$ and select time $(t',j')$ before this jump and show item (i) of Definition~\ref{defstabgraf} holds.
\\
First, observe that $z_1\in \widehat D$ holds also in the case where $D$ is unbounded following \eqref{eqphis} and (iv). Hence, $(z_1,z_2)\in \mathcal{A}_{01}$ holds.
To prove $\tilde{\jmath}>0$ in \eqref{eqphis}, suppose the contrary, i.e.\ $\tilde{\jmath}=0$.
Let $t^\star\leq t$ denote the minimum continuous time such that $\inf_{z\in \widehat D}\|(\phi^\star(\tau,j)-z,\phitwo(\tau,0)-G(z))\|<s$ for $\tau\in (t^\star,t]$ and $[t^\star,j]\in \dom \phi^\star$. Continuity of hybrid arcs during flow either implies $(t^\star,j-1)\in \dom \phi^\star$ or there exists $z_1^\star \in \widehat D$ such that
$\|(\phi^\star(t^\star,j)-z_1,\phitwo(t^\star,0)-G(z_1^\star))\|= s.
$
The first option implies $\phi(t^\star,j)\in G(\widehat D)$ such that $\|\phi^\star(t^\star,j)-z\|\geq \min_{u\in \widehat D,w\in G(\widehat D)}\|u-w\|>s$, contradicting $\inf_{z\in \widehat D}\|(\phi^\star(\tau,j)-z,\phitwo(\tau,0)-G(z))\|<s$. Otherwise, by design of $\bar s$, we find
\beq{
\inf_{z\in C\cup D}\|(\phi^\star(t^\star,j)-z,\phitwo(t^\star,0)-z)\|> \bar s>s\label{eqpfth2bars}
}
and $\inf_{z\in D}\|(\phi^\star(t^\star,j)-G(z),\phitwo(t^\star,0)-z)\|> \bar s>s$. Again exploiting continuity of hybrid arcs during flow, there cannot exist a hybrid time interval $[\tau^\star,t^\star)\times \{j\}$, with $\tau^\star<t^\star$ and $\inf_{z\in C\cup D}\|(\phi^\star(\tau,j)-z,\phitwo(\tau,0)-z)\|<s$ or $\inf_{z\in D}\|(\phi^\star(\tau,j)-G(z),\phitwo(\tau,0)-z)\|<s$ for $\tau\in [\tau^\star,t^\star)$ and, since (vi) holds, the only remaining option is $(t^\star,j)=(0,0)$, in which (v) contradicts \eqref{eqpfth2bars}. A contradiction is found in every scenario and  $\tilde{\jmath}>0$.
\\
Since \eqref{eqphis} implies $\|\phitwo(t,\tilde{\jmath})-G(z_1)\|\leq s$, $s<K$ and $z_1\in \widehat D$ has been obtained above, we find $\phitwo(t,\tilde{\jmath})\in G(D\cap \mathcal{B}_{2K})+\mathcal B_{K}$, such that $\|\phitwo(t,\tilde{\jmath})\|<\bar K$ follows from the construction of $\bar K$.
As, in addition, $\phitwo(t,\tilde{\jmath})\in G(D)+\mathcal{B}_s$ holds, $s<\epsilon_3$ and $\epsilon_3$ is selected as in Lemma~\ref{lempostj} with $\epsilon_4=\tfrac{\min(\varepsilon,2r)}{2(1+\tilde F)}$, there exists $t'\in [t-\tfrac{\min(\varepsilon,2r)}{2(1+\tilde F)},t]$ such that $\phitwo(t',\tilde{\jmath})\in G(D)$. Similarly, we infer that $\inf_{z\in \widehat D}\|(\phi^\star(\tau,j)-z,\phitwo(\tau,\tilde{\jmath})-G(z))\|\leq s$ for $\tau\in [t^\star,t]$ and $t^\star=\max(t',\min\{t\in \mathbb{R}\,:\,  (t,j)\in \dom \phi^\star\})$. Hence,
\beq{
\phi^\star(t^\star,j)\in D+\mathcal{B}_s\label{eqphi1tprime}
} is obtained, which implies $\phi^\star(t^\star,j)\notin G(D)$ and $t^\star=t'$.
\\
From $\phitwo(t',\tilde{\jmath})\in G(D)$, $\tilde{\jmath}\geq 1$ and item (iii), we find $(t',\tilde{\jmath}-1)\in \dom \phitwo$. Since $|t-t'|\leq \tfrac{\min(\varepsilon,2r)}{2(1+\tilde F)}<\varepsilon$, we will conclude this case and show that $(t',j')$, with $j'=\tilde{\jmath}-1$, satisfies $\|\phi^\star(t,j)-\phitwo(t',j')\|<\varepsilon$.
We first prove
\beq{
\rho_{\mathcal{A}_{00}}(\phi^\star(t',j),\phitwo(t',\tilde{\jmath}-1))\leq \tfrac{\varepsilon}2
\label{eq-pfthm-tprimejprimeminus}
} holds by considering the case of $t'=0$ separately, followed by the case in which $t'>0$. If $t'=0$, we use $\phi^\star(t',j)\notin G(D)$ obtained above to deduce $j=0$ and $\phitwo(t',\tilde{\jmath}-1)\in D$ to deduce $\tilde{\jmath}-1=0$ (since $G(D)\cap D=\emptyset$ by item~(i)), such that $\|\phi^\star(t',j)-\phitwo(t',\tilde{\jmath}-1)\|=\|\phi^\star(0,0)-\phitwo(0,0)\|\leq s$ by (v). As $\rho_{\mathcal{A}_{00}}(x,y)\leq \|x-y\|$ for all $x,y\in C\cup D$ and $s\leq \tfrac{\varepsilon}2$, we obtain \eqref{eq-pfthm-tprimejprimeminus}.
If $t'>0$, using items (i) and (iii) and the inclusion $(t',\tilde{\jmath}-1)\in \dom \phitwo$, we observe that there exists a time $t''<t'$ such that for $\tau\in (t'',t'),$ the equality
 \begin{equation}\label{eq-pfthm-setjs}
\{\jmath\in \{0,1,\ldots\} \,:\,   (\tau,\jmath)\in \dom \phitwo\}=\{\tilde{\jmath}-1\}
\end{equation}
holds, i.e., no jumps of $\phitwo$ occur in the open continuous-time interval $(t'',t')$.
From \eqref{eqphi1tprime}
we find that  $(\tau,j)\in \dom \phi^\star$ holds for all $\tau\in (t''',t'),$ and some $t'''\in [t'',t')$.
For $\tau\in (t''',t')$, (vi) implies $\rho_\mathcal{A}(\phi^\star(\tau,j),\phitwo(\tau,\tilde{\jmath}-1))\!<\!s$ for $\tau\!\in\! (t''',t')$. Hence, for a sequence $\{\tau_k\}_{k\in \mathbb{N}}$ with $\tau_k\!\in\! (t''',t')$ and $\lim_{k\to \infty}\!\tau_k\!=\!t'$ we find $\lim_{k\to \infty }\phi^\star(\tau_k,j)\in D+\mathcal{B}_s$, $\lim_{k\to \infty}\phitwo(\tau_k,\tilde{\jmath}-1)\in D$ and $\lim_{k\to\infty}\rho_\mathcal{A}(\phi^\star(\tau_k,j)$, $\phitwo(\tau_k,\tilde{\jmath}-1))\!\leq\! s$. For each sufficiently large $k$, $\rho_\mathcal{A}(\phi^\star(\tau_k,j),\phitwo(\tau_k,\tilde{\jmath}-1))=\rho_{\mathcal{A}_{00}}(\phi^\star(\tau_k,j),\phitwo(\tau_k,\tilde{\jmath}-1))$. Namely, if $x_1\!\in\! D\!+\!\mathcal{B}_s$ and $x_2\!\in\! D$, then $\rho_{\mathcal{A}}(x_1,x_2)\leq s$ implies $\rho_{\mathcal{A}}(x_1,x_2)=\rho_{\mathcal{A}_{00}}(x_1,x_2)$.
 By continuity of $\rho_{\mathcal{A}_{00}}$ and continuity of the hybrid arcs for fixed $j$, $\tilde{\jmath}$, we find
$\lim_{k\to \infty}\rho_{\mathcal{A}_{00}}(\phi^\star(\tau_k,j),\phitwo(\tau_k,\tilde{\jmath}-1))=\rho_{\mathcal{A}_{00}}(\phi^\star(t',j),$ $ \phitwo(t',\tilde{\jmath}-1))\leq s$. With $\rho_{\mathcal{A}_{00}}(\phi^\star(t',j),$ $ \phitwo(t',\tilde{\jmath}-1))\geq \tfrac{1}{\sqrt 2 } \|\phi^\star(t',j)-\phitwo(t',\tilde{\jmath}-1)\|$ and $s<\tfrac \varepsilon {2\sqrt 2\tilde F}$, we find \eqref{eq-pfthm-tprimejprimeminus}.
\\
From $\rho_{\mathcal{A}_{01}}(\phi^\star(t,j),\phitwo(t,\tilde{\jmath}))=\inf_{z\in D}\|(\phi^\star(t,j)-z,\phitwo(t,\tilde{\jmath})-G(z))\|< s$, we find $\phi^\star(t,j)\in D+\mathcal{B}_s$ and, since $s\leq \tfrac{r}{\tilde F+1}$ and $|t'-t|\leq \tfrac{\min(\varepsilon,2r)}{2(1+\tilde F)}\leq \tfrac{r}{\tilde F+1}$, we obtain $\dot \phi^\star(\tau,j)\leq \tilde F$ for $\tau\in [t',t]$.
Exploiting
$|t'-t|\leq \tfrac{\varepsilon}{2(1+\tilde F)}$
and $\tfrac{\tilde F}{1+\tilde F}<1$, we get $\|\phi^\star(t',j)-\phi^\star(t,j)\|<|t'-t|\tilde F\leq \tfrac{\varepsilon}{2}$. With \eqref{eq-pfthm-tprimejprimeminus} and $j'=\tilde{\jmath}-1$, we find   $\|\phi^\star(t,j)-\phitwo(t',j')\|\leq  \|\phi^\star(t,j)-\phi^\star(t',j)\|+
\|\phi^\star(t',j)-\phitwo(t',j')\|<\varepsilon$,
such that $(t',j')$ satisfy the theorem conditions.
\\
\textbf{Case 3: $z_2\in D,z_1= G(z_2)$.\quad }
Since $\phitwo(t,\tilde \jmath)$ is close to $D$, we will apply Lemma~\ref{lemprej} to prove a jump of $\phitwo$ will occur soon, and select $(t',j')$ directly after this jump. Subsequently, we conclude this case by showing that item (i) of Definition~\ref{defstabgraf} holds for this hybrid time instant.
\\
For this purpose, first, we observe that $z_2\in \widehat D$ holds also in the case where $D$ is unbounded. Namely, as \eqref{eqphis} implies $\|\phi^\star(s)-G(z_2)\|<s$ and $s<K$, we find with (iv) that $\|G(z_2)\|<2K$. Hence, $z_2\in G^{-1}(D\cap \mathcal{B}_{2K})\subseteq \mathcal{B}_{\bar K}$ holds by construction of $\bar K$ and $(z_1,z_2)\in \mathcal{A}_{10}$ is verified.
Since $\|\phitwo(t,\tilde{\jmath})-z_2\|<s<K$ follows from \eqref{eqphis}, we find $\|\phitwo(t,\tilde{\jmath})\|<\bar K$ and $\phitwo(t,\tilde{\jmath})\in D+\mathcal{B}_s$. Since $s<\epsilon_1$, we can apply Lemma~\ref{lemprej} and conclude there exists a time $t'\in [t,t+\tfrac{\min(\varepsilon,2r)}{2(1+\tilde F)}]$ such that $\phitwo(t',\tilde{\jmath})\in \widehat D$. With items (i) and (ii), we find $(\tau,\tilde{\jmath}+1)\in \dom \phitwo$ for $\tau\in [t',t'')$, with some $t''>t'$. Reasoning analogously as in Case 2, we obtain $(t',j)\in \dom \phi^\star$, $\phi^\star(t',j)\in G(D)+\mathcal{B}_s$, such that $\phi^\star(t',j)\notin D$ and, choosing $t'''$ sufficiently small, we find $(\tau,j)\in \dom \phi^\star$ for $\tau\in [t',t''')$. Taking a sequence $\{\tau_k\}_{k\in \mathbb{N}}$ with $\tau_k>t'$ and $\lim_{k\to \infty} \tau_k=t'$,  we find $\lim_{k\to \infty}\rho_{\mathcal{A}_{00}}(\phi^\star(\tau_k,j),\phitwo(\tau_k,\tilde{\jmath}+1))=\rho_{\mathcal{A}_{00}}(\phi^\star(t',j),\phitwo(t',\tilde{\jmath}+1))\leq s$. As $s<\tfrac \varepsilon {2\sqrt 2\tilde F}$ and $\tilde F\geq 1$, we find $\|\phi^\star(t',j)-\phitwo(t',\tilde{\jmath}+1)\|<\tfrac {\varepsilon}{2}$.
\\
 From \eqref{eqphis} and $z_2\in \widehat D$, we find $\phi^\star(t,j)\in G(\widehat D)+\mathcal{B}_s$ and, since $s\!<\!\tfrac{r}{\tilde F+1}$, we obtain $\dot \phi^\star(\tau,j)\!\leq\! \tilde F$ for $\tau\!\in\! [t,t']$, since $|t'-t|<\tfrac{\min(\varepsilon,2r)}{2(1+\tilde F)}<\tfrac{r}{\tilde F+1}$. We deduce $\|\phi^\star(t',j)-\phi^\star(t,j)\|<|t-t'|\tilde F\leq \tfrac{\varepsilon}2$ from $|t-t'|\leq
\tfrac{\varepsilon}{2(1+\tilde F)}$. Selecting $j'=\tilde{\jmath}+1$, we obtain $\|\phi^\star(t,j)-\phitwo(t',j')\|\leq
\|\phi^\star(t,j)-\phitwo(t',j)\|+\|\phi^\star(t',j)-\phitwo(t',j')\|\leq  \varepsilon$, and the conclusion of the theorem is verified.
\\
As $(t',j')$ has been constructed for all three cases and arbitrary $(t,j)\in \dom \phi^\star$, the theorem is proven.
\hfill $\Box$
\end{pf}

\def\Elproofname{PROOF of THEOREM~\ref{thm-dafis}.}
\begin{pf}
We exploit items (i) and (ii) of Definition~\ref{dfnstabtrajrhoa} and Theorem~\ref{thm-mismatch} to conclude stability in graphical sense as defined in Definition~\ref{defstabgraf}, and exploit the combination of (iii), (iv) of Definition~\ref{dfnstabtrajrhoa} and Theorem~\ref{thm-mismatch} to conclude (asymptotic) stability in graphical sense. Given a pair of solutions $(\phi^\star,\phitwo)$, Theorem~\ref{thm-mismatch} only provides statements for all $(t,j)\in \dom \phi^\star$, (see item (i) of Definition~\ref{defstabgraf}). Statement (ii) of Definition~\ref{defstabgraf} is attained by another application of Theorem~\ref{thm-mismatch} for the solution pair $(\phi^{{\star'}},\phitwo')$, which we select as $(\phitwo,\phi^\star)$. (Asymptotic) stability of $\phi^\star$ with respect to $\rho_{\mathcal A}$
will be used to show that (vi) of Theorem~\ref{thm-mismatch} holds for both solution pairs $(\phi^\star,\phitwo)$ and $(\phi^{\star'},\phitwo')$. 
\\
Consider system \eqref{eqhs} and solution $\phi^\star$ satisfying the conditions of Theorem~\ref{thm-dafis}. We select $K>0$ such that either $\|\phi^\star(t,j)\|<K$ holds for all $(t,j)\in \dom \phi^\star$ or $\|x\|<K$ for $x\in D$, cf. (iv) of Theorem~\ref{thm-mismatch}. Exploiting also item (i) in Theorem~\ref{thm-mismatch} and local boundedness of $G$ if $D$ is unbounded, we can select a scalar $K'$ such that the combination of the inequality  $\rho_{\mathcal{A}}(\phi^\star(t,j),y)<K$ for any $y\in D$ and item (iv) of Theorem~\ref{thm-mismatch} implies $\|y\|<K'$ for any $(t,j)\in \dom \phi^\star$.
Let $\phi^\star$ be stable with respect to $\rho_{\mathcal A}$. To prove existence of a scalar $\delta>0$ for any $\varepsilon>0$ (see Definition~\ref{defstabgraf}) such that (i) and (ii) in Definition~\ref{defstabgraf} hold, we first fix an arbitrary $\varepsilon>0$. By application of Theorem~\ref{thm-mismatch}, we find a scalar $s>0$ such that for any solution $\phitwo$ that satisfies (v)-(vi) of Theorem~\ref{thm-mismatch}, item (i) of Definition~\ref{defstabgraf} holds.
Considering the pair $(\phi^{\star'},\phitwo')=(\phitwo,\phi^\star)$ of solutions, we apply Theorem~\ref{thm-mismatch} and find $s'>0$ such that if $\phi$ satisfies (v)-(vi), with $s$ replaced by $s'$,
for every $(t',j')\in \dom \phitwo,$ there exists a hybrid time $(t,j)\in \dom \phi^\star$, with $|t-t'|<\varepsilon$, such that $\| \phi^\star(t,j)-\phitwo(t',j')\|<\varepsilon$.
We now select $\delta=\delta_w(\min(s,s',K))$, with $\delta_w(\cdot)$ given in Definition~\ref{dfnstabtrajrhoa}, and consider an arbitary solution $\phitwo$ with $\|\phi^\star(0,0)-\phitwo(0,0)\|<\delta$ (cf.\  Definition~\ref{defstabgraf}) and will show that conditions (i) and (ii) in Definition~\ref{defstabgraf} hold. We note that  $\rho_{\mathcal{A}}(\phi^\star(0,0),\phitwo(0,0))\leq \|\phi^\star(0,0)-\phitwo(0,0)\|<\delta\leq \delta_w(s)$ implies that item (i) of Definition~\ref{dfnstabtrajrhoa} is verified with $\varepsilon_w=s$. Hence, for the solution pair $(\phi^\star,\phitwo)$ with scalar $s$, condition (v) of Theorem~\ref{thm-mismatch} holds and, by item (i) of Definition~\ref{dfnstabtrajrhoa}, we conclude that (vi) of Theorem~\ref{thm-mismatch} holds. As (iv) holds by assumption, we apply Theorem~\ref{thm-mismatch} and conclude item (i) of Definition~\ref{defstabgraf}.
\\
We note that  $\rho_{\mathcal{A}}({\phi^\star}(0,0),\phitwo(0,0))\leq \|\phi^\star(0,0)-\phitwo(0,0)\|<\delta_w(s')$ implies that item (ii) of Definition~\ref{dfnstabtrajrhoa} is verified with $\varepsilon_w=s'$. Hence, aiming to apply Theorem~\ref{thm-mismatch} with the solution pair $({\phi^\star}',\phitwo')=(\phitwo,\phi^\star)$ and scalar $s'$, we observe that condition (v) holds and, by item (ii) of Definition~\ref{dfnstabtrajrhoa}, we conclude that (vi) of Theorem~\ref{thm-mismatch} holds (also for the pair  $({\phi^\star}',\phitwo')$). As (iv) of the same theorem holds by assumption, we can apply Theorem~\ref{thm-mismatch} to conclude that item (i) of Definition~\ref{defstabgraf} holds for the solution pair $({\phi^\star}',\phitwo')$. As a direct consequence, item (ii) of Definition~\ref{defstabgraf} holds for the solution pair $(\phi^\star,\phitwo)$. Since we have proven (i) and (ii) for this solution pair, and $\phitwo$ is selected arbitrarily, the solution $\phi^\star$ is stable in graphical sense.
\\
We now show \emph{asymptotic stability}. Assume that $\phi^\star$ is asymptotically stable with respect to $\rho_{\mathcal{A}}$, let $r=r_w>0$ be as in Definition~\ref{dfnstabtrajrhoa}, and select $\varepsilon>0$ arbitrarily. Given $\varepsilon$, let $s$ be as in Theorem~\ref{thm-mismatch} and let $s'$ be selected as above. In addition, consider any $\phitwo$ with $\|\phi^\star(0,0)-\phitwo(0,0)\|\leq r$.
\\
Applying Lemma~\ref{lemprej} with $\epsilon_2=\varepsilon$, we find some positive scalar $\epsilon_1$. Let $\bar s>0$ be as in the proof of Theorem~\ref{thm-mismatch}.
We now consider Definition~\ref{dfnstabtrajrhoa} with $\varepsilon_w=\min(\tfrac{s}{\sqrt 2},\tfrac{s'}{\sqrt 2},\epsilon_1,\bar s)$, and find that there exists a time $T_w >0$ such that items (iii) and (iv) of Definition~\ref{dfnstabtrajrhoa} hold. In particular, this implies that there exist $J_w,J_w'$ such that $\rho_{\mathcal{A}}(\phi^\star(T_w,J_w),\phitwo(T_w,J'_w))<\varepsilon_w$. With $\varepsilon_w\leq \epsilon_1$ and Lemma~\ref{lemprej}, we conclude that there exist hybrid times $(T,J)\in \dom \phi^\star$ and $(T,J')\in \dom \phitwo$, with $T\in [T_w,T_w+\varepsilon]$, such that $\rho_{\mathcal{A}}(\phi^\star(T,J),\phitwo(T,J'))=\rho_{\mathcal{A}_{00}}(\phi^\star(T,J),\phitwo(T,J'))$, where the design of $\bar s$ and $\varepsilon_w\leq \bar s$ are used. From $\rho_{\mathcal{A}_{00}}(\phi^\star(T,J),\phitwo(T,J')<\varepsilon_w$, we conclude $\|\phi^\star(T,J)-\phitwo(T,J')\|<\min(s,s')$.
\\
Let $(\phi_s^\star,\phitwo_s)$ be constructed such that $\phi^\star(T+t,J+j)=\phi_s(t,j)$ for all $(t,j)\in \dom \phi^\star_s$ and $\phitwo(T+t,J'+j)=\phitwo_s(t,j')$ for all $(t,j')\in \dom \phitwo_s$. For $(\phi^\star_s,\phitwo_s)$ and $(\phitwo_s,\phi^\star_s)$, all conditions of Theorem~\ref{thm-mismatch} hold, such that items (i),(ii) of Definition~\ref{defstabgraf} follow for $(\phi^\star_s,\phitwo_s)$ and $(\phitwo_s,\phi^\star_s)$. We conclude that the scalar $T>0$ constructed above ensures items (iii) and (iv) of Definition~\ref{defstabgraf}. Consequently, $\phi^\star$ is \emph{asymptotically stable in graphical sense}.\hfill $\Box$
\end{pf}
\end{document}